\documentclass[a4paper,11pt]{article}
\usepackage[margin=1.1in]{geometry}
\usepackage[T1]{fontenc}
\usepackage{lmodern}
\usepackage[expansion=false]{microtype}
\usepackage{amsmath,amssymb,amsthm}
\usepackage[section]{placeins}
\usepackage{needspace}
\usepackage{tikz}
\usetikzlibrary{arrows.meta,calc}
\usepackage[hidelinks]{hyperref}
\newtheorem{theorem}{Theorem}
\newtheorem{proposition}[theorem]{Proposition}
\newtheorem{lemma}[theorem]{Lemma}
\hypersetup{pdftitle={A linear edge bound for k-quasi-planar graphs},
  pdfauthor={Ankan Sadhu}}
\tikzset{
  diagram/.style={font=\small,line cap=round,line join=round},
  first curve/.style={black,line width=.95pt},
  second curve/.style={black,line width=.85pt,densely dashed},
  surface border/.style={black!65,line width=.45pt},
  map arrow/.style={-{Latex[length=2mm]},line width=.55pt}
}

\title{A linear edge bound for $k$-quasi-planar graphs}

\author{Ankan Sadhu\\[.4em]
  \small Department of Computer Science and Engineering\\
  \small Government College of Engineering and Ceramic Technology\\
  \small Kolkata, India\\[.25em]
  \small\href{mailto:ankan-sadhu-br24-3003@gcect.ac.in}
    {\texttt{ankan-sadhu-br24-3003@gcect.ac.in}}}
\date{}

\begin{document}

\maketitle

\begin{center}
\small\textbf{Supersession notice.} This preprint has been superseded by the joint paper
\emph{Polynomial bounds for quasiplanar graphs},
\href{https://arxiv.org/abs/2611.21416}{arXiv:2611.21416}.
\end{center}

\begin{abstract}
We prove the $k$-quasi-planar graph conjecture for topological
drawings. For every integer $k\ge 2$, a simple graph on $n\ge 3$
vertices admitting a drawing with no $k$ pairwise crossing edges
has at most $1056(k-1)^2(n-2)$ edges. Thus the number of edges
is linear in $n$, with an explicit quadratic dependence on $k$.
The result allows multiple crossings between the same pair of
edges, including edges with a common endpoint. The proof uses
the double of a sphere with boundary to represent graph edges
by closed geodesics on a hyperbolic surface. The endpoints of
geodesic lifts induce partial orders of width at most $k-1$.
A stability bound for greedy chain partitions, together with
a counting argument over a triangulation, gives the stated
edge bound.
\end{abstract}

\begin{quote}

\end{quote}

\section{Introduction}

The quasi-planar graph conjecture asks whether, for every fixed
integer $k$, a simple graph drawn without $k$ pairwise crossing edges
has only linearly many edges. Pach, Shahrokhi, and
Szegedy~\cite{PSS96} studied the problem for straight-line drawings.
Here we prove the conjecture for topological drawings, with an
explicit coefficient quadratic in $k$.

All graphs in this paper are finite and simple. In a \emph{topological
drawing}, vertices are represented by distinct points in the plane,
and edges by Jordan arcs whose interiors avoid all vertices. Every
intersection of two edge interiors is a \emph{proper crossing}: the
two branches of one edge alternate with the two branches of the other.
We require finitely many crossings and no point where three edge
interiors meet. For an integer $k\ge 2$, the drawing is
\emph{$k$-quasi-planar} if no $k$ distinct edges cross pairwise in
their interiors. This definition allows edges to cross repeatedly,
including edges that share an endpoint.

\Needspace{7\baselineskip}
\begin{theorem}\label{thm:main}
Let $k\ge 2$ and $n\ge 3$ be integers. Every simple graph on
$n$ vertices that admits a $k$-quasi-planar topological drawing
has at most
\[
 1056(k-1)^2(n-2)
\]
edges. In particular, for each fixed $k$, the number of edges
is at most $C_k n$, where $C_k=1056(k-1)^2$.
\end{theorem}

For $k=2$, the conjecture reduces to the planar bound of $3n-6$
edges for $n\ge 3$. Capoyleas and
Pach~\cite{CapoyleasPach1992} proved a linear bound for straight-line
drawings whose vertices are in convex position. Pach,
Radoi\v{c}i\'c, and T\'oth~\cite{PachRadoicicToth2006} established
linearity for $k=3$, extending the result of Agarwal and his
coauthors~\cite{AgarwalEtAl1997} for drawings in which each pair
of edges meets at most once. Ackerman and
Tardos~\cite{AckermanTardos2007} improved the bound for $k=3$,
and Ackerman~\cite{Ackerman2009} established linearity for $k=4$.
For general $k$, Fox, Pach, and Suk~\cite{FoxPachSuk2024} proved
an upper bound of $n(\log n)^{O(\log k)}$ edges.

\paragraph{Related preprints and the present bound.}
Version 1 of this manuscript~\cite{Sadhu26}, posted on
1 November 2026, derived fixed-$k$ linearity from the
surface-pattern theorem of Beeker and Lazarovich~\cite{BL17}.
Pach, Tardos, and Zeng~\cite{PTZ26} posted an independent
combinatorial proof on 5 November 2026. They report encountering
Version 1 while preparing their manuscript and obtain the bound
$2^{(2+o(1))k\log_2 k}n$ as $k\to\infty$. The present version
gives the explicit bound in Theorem~\ref{thm:main} by a different
argument: a hyperbolic-geometric construction followed by a
quantitative stability estimate for greedy partitions into chains.

Pach, Tardos, and Zeng~\cite{PTZ26} conjecture a bound of
$O(kn)$ edges with an absolute implied constant. When $2k-1$
divides $n$, disjoint convex straight-line drawings of copies of
$K_{2k-1}$ have $(k-1)n$ edges. Such drawings are
$k$-quasi-planar: $k$ pairwise crossing chords require $2k$
distinct endpoints.

We begin by removing disks around the vertices, as in the proof
of Theorem~1.7 of Fox, Pach, and
Suk~\cite{FoxPachSuk2024}. The remaining portions of the edges
are arcs on a sphere with boundary. Doubling the surface and
straightening the doubled arcs produces closed geodesics. On each
lifted triangulation side, the ideal endpoints of the geodesics
determine two linear orders. Their common partial order has width
at most $k-1$, because incomparable geodesics intersect. The
orders agree on geodesics that cross two sides directed both toward,
or both away from, their common vertex.

For the final count, we use the triangulation-and-counting
framework of Pach, Tardos, and Zeng~\cite[Sections~2 and~4]{PTZ26}.
Starting at a boundary circle, we follow the arcs of two
consecutive elements in a greedy chain. The two graph edges have
different other endpoints, so their arcs must eventually separate.
A potential function bounds the number of consecutive pairs lost
when the arcs pass between triangulation sides. We thereby assign
every graph edge to a record in one of $8(n-2)$ triangles, and
bound the number of records per triangle by $O(k^2)$.

\section{From a drawing to closed curves}\label{sec:double}

We first turn the drawing edges into closed curves. Intersection
parity will distinguish curves belonging to different edges.
A closed curve on a surface is \emph{essential} if it is not
homotopic to a point. A \emph{free homotopy} does not fix a
basepoint. An \emph{unoriented free homotopy class} identifies a
closed curve with its reverse and with all curves freely
homotopic to it.

\begin{proposition}\label{prop:double}
Let $G$ be a simple graph with $n\ge 3$ vertices, equipped
with a topological drawing in the plane. There exist a closed
orientable surface $S$ of genus $n-1$ and an essential simple
closed curve $\gamma_e$ on $S$ for each edge $e\in E(G)$.
The curves represent distinct unoriented free homotopy classes.
Two such curves intersect if and only if the corresponding
edges cross in the drawing.
\end{proposition}

\begin{proof}
Add the point at infinity to the plane, obtaining a sphere.
Treat each crossing as an additional vertex; the drawing then
gives a finite graph embedded in the sphere. Around each original
vertex $v$, choose a small closed disk $D_v$. Choose the disks to
be pairwise disjoint, to avoid all crossings, and to meet each
incident edge only in an initial subarc. If $v$ is isolated,
choose $D_v$ to meet the drawing only at $v$.

Remove the interiors of all the disks, and call the resulting
surface $P$. Write $b_v=\partial D_v$ for its boundary circle
at $v$. An arc in $P$ is \emph{properly embedded} if it is
simple, has both endpoints on the boundary, and has its
interior in the interior of $P$. For an edge $e$ joining vertices
$u$ and $v$, let $a_e$ be the portion of the edge in $P$.
This arc is properly embedded and joins $b_u$ to $b_v$.
All arc endpoints are distinct, including those from edges
incident with the same vertex. The arc interiors retain
exactly the original crossings (Figure~\ref{fig:truncate}).

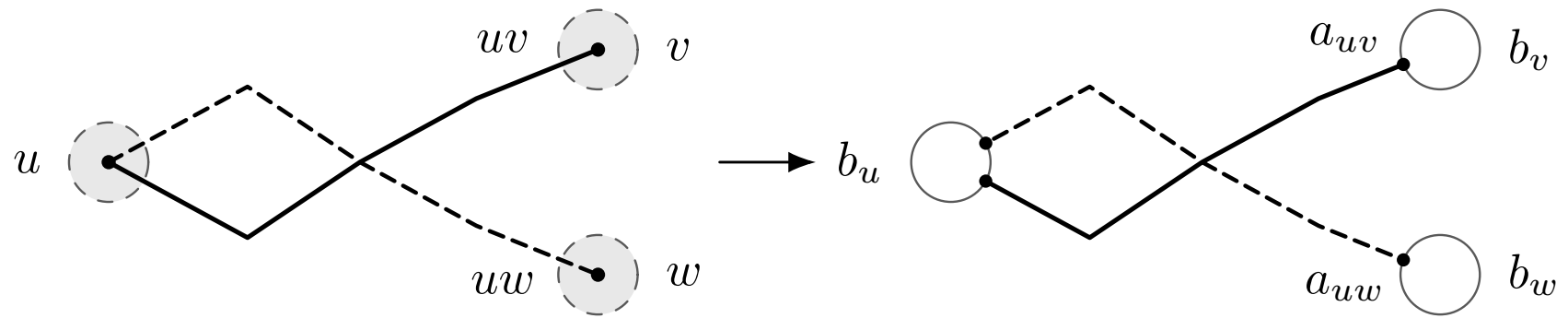
\begin{figure}[htbp]
\centering
\begin{tikzpicture}[diagram,x=.95cm,y=.95cm]
  \foreach \panel/\offset in {0/0,1/6.05} {
    \begin{scope}[xshift=\offset cm]
      \coordinate (u) at (0,1);
      \coordinate (v) at (3.7,1.85);
      \coordinate (w) at (3.7,.15);
      \coordinate (a) at (1.05,.43);
      \coordinate (b) at (1.05,1.57);
      \coordinate (c) at (2.78,1.48);
      \coordinate (d) at (2.78,.52);
      \coordinate (z) at (1.9,1);
      \ifnum\panel=0
        \foreach \p in {u,v,w}
          \fill[black!9] (\p) circle[radius=.30];
      \fi
      \draw[first curve] (u) -- (a) -- (z) -- (c) -- (v);
      \draw[second curve] (u) -- (b) -- (z) -- (d) -- (w);
      \ifnum\panel=0
        \foreach \p in {u,v,w} {
          \draw[surface border,dashed] (\p) circle[radius=.30];
          \fill (\p) circle[radius=1.45pt];
        }
        \node[left=10pt] at (u) {$u$};
        \node[right=10pt] at (v) {$v$};
        \node[right=10pt] at (w) {$w$};
        \node[above] at (2.98,1.69) {$uv$};
        \node[below] at (2.98,.31) {$uw$};
        \node at (1.85,-.48) {(a) Vertex disks};
      \else
        \foreach \p in {u,v,w}
          \draw[surface border,fill=white] (\p) circle[radius=.30];
        \foreach \p/\q in {u/a,u/b,v/c,w/d}
          \fill ($(\p)!.285cm!(\q)$) circle[radius=1.35pt];
        \node[left=10pt] at (u) {$b_u$};
        \node[right=10pt] at (v) {$b_v$};
        \node[right=10pt] at (w) {$b_w$};
        \node[above] at (2.98,1.69) {$a_{uv}$};
        \node[below] at (2.98,.31) {$a_{uw}$};
        \node at (1.85,-.48) {(b) Truncated edges};
      \fi
    \end{scope}
  }
  \draw[map arrow] (4.62,1) -- (5.35,1);
\end{tikzpicture}
\caption{The solid and dashed curves represent $uv$ and $uw$.
After the vertex disks are removed, the arcs have distinct endpoints
on $b_u$ and still cross in the interior of~$P$.
\emph{TikZ figure prepared with OpenAI Codex assistance.}}
\label{fig:truncate}
\end{figure}

Take two copies $P_+$ and $P_-$ of $P$ with opposite orientations.
Identify corresponding points on their boundaries. The resulting
closed, connected, orientable surface $S$ is the \emph{double} of $P$.
For each vertex $v$, let $s_v$ be the image of $b_v$ in $S$.
Write $\chi(X)$ for the Euler characteristic of a surface $X$.
Since $P$ is a sphere with $n$ open disks removed, $\chi(P)=2-n$.
Every gluing circle has Euler characteristic zero, and hence
\[
  \chi(S)=2\chi(P)=4-2n,
\]
and $S$ has genus $n-1$.

For each edge~$e$, let $a_e^+$ and $a_e^-$ be the copies of $a_e$
in $P_+$ and $P_-$. Their union in $S$ is a simple closed
curve~$\gamma_e$ (Figure~\ref{fig:glue}). Distinct arcs have
distinct boundary endpoints, so doubled curves meet only at copies
of crossings in the drawing.

\begin{figure}[htbp]
\centering
\begin{tikzpicture}[diagram]
  \path[fill=black!10] (2.55,0)
    arc[start angle=0,end angle=180,x radius=2.55,y radius=1.35]
    -- (-1.38,0)
    arc[start angle=180,end angle=0,x radius=1.38,y radius=.35] -- cycle;
  \path[fill=black!3] (-2.55,0)
    arc[start angle=180,end angle=360,x radius=2.55,y radius=1.35]
    -- (1.38,0)
    arc[start angle=360,end angle=180,x radius=1.38,y radius=.35] -- cycle;
  \draw[surface border] (0,0) ellipse[x radius=2.55,y radius=1.35];
  \draw[surface border] (0,0) ellipse[x radius=1.38,y radius=.35];
  \draw[black!70,line width=.6pt,densely dotted] (-2.55,0) -- (-1.38,0);
  \draw[black!70,line width=.6pt,densely dotted] (1.38,0) -- (2.55,0);
  \draw[first curve] (0,0) ellipse[x radius=1.95,y radius=.85];
  \fill (-1.95,0) circle[radius=1.6pt] (1.95,0) circle[radius=1.6pt];
  \node[left=4pt] at (-2.55,0) {$s_u$};
  \node[right=4pt] at (2.55,0) {$s_v$};
  \node at (0,1.10) {$P_+$};
  \node at (0,-1.10) {$P_-$};
  \node at (0,.60) {$a_e^+$};
  \node at (0,-.60) {$a_e^-$};
\end{tikzpicture}
\caption{An annular neighborhood of the doubled curve $\gamma_e$
for $e=uv$. The copies $a_e^+$ and $a_e^-$ join at the marked points.
The dotted segments are the portions of $s_u$ and $s_v$ in the
annulus.
\emph{TikZ figure prepared with OpenAI Codex assistance.}}
\label{fig:glue}
\end{figure}
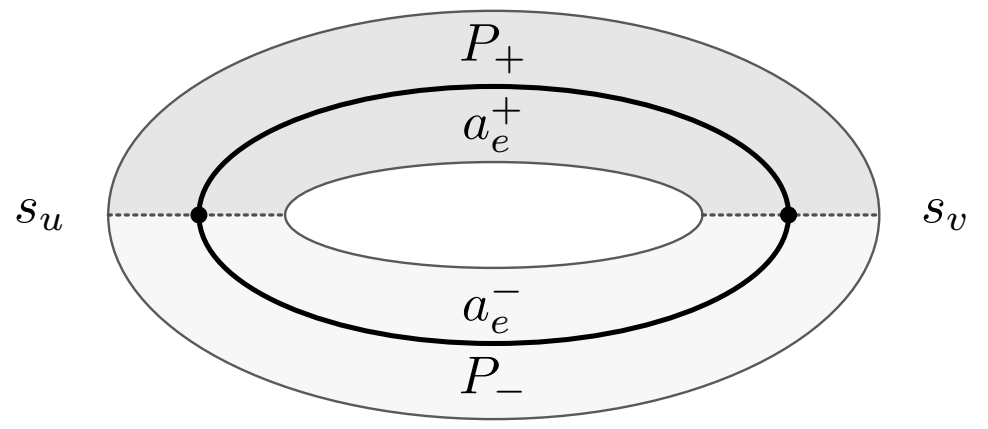

For transverse closed curves $\gamma$ and $s$ on $S$, their
\emph{mod-two intersection number} $I_2(\gamma,s)$ is the parity
of their number of intersection points. A generic homotopy creates
or cancels intersections in pairs. Thus $I_2(\gamma,s)$ depends
only on the free homotopy classes of the curves. In particular,
a contractible curve has intersection number zero: it can be
contracted to a point outside the other curve.

If $e$ has endpoints $u$ and $v$, then $\gamma_e$ crosses $s_u$
and $s_v$ once each and misses every other gluing circle. Thus
\[
  I_2(\gamma_e,s_x)=
  \begin{cases}
    1 & \text{if } x\in\{u,v\},\\
    0 & \text{otherwise}.
  \end{cases}
\]
Since $I_2(\gamma_e,s_u)=1$, the curve $\gamma_e$ is essential.
Distinct edges of a simple graph have different unordered endpoint
pairs. For any two such edges, some gluing circle therefore
distinguishes the mod-two intersection numbers of their curves.
These numbers are unchanged by reversing a curve. Hence the
curves represent distinct unoriented free homotopy classes.
\end{proof}

\section{Geodesics and their endpoint orders}\label{sec:orders}

The universal cover of a closed hyperbolic surface is the
hyperbolic plane. In the Klein disk model, complete geodesics
are open chords. The two points where a chord meets the boundary
circle are its \emph{ideal endpoints}. Two geodesics with four
distinct ideal endpoints intersect exactly when their endpoints
alternate around the circle. A \emph{lift} of an embedded closed
curve is a component of its preimage in the universal cover.

Throughout the proof, a \emph{triangulation} means a decomposition
into triangular cells; distinct sides may share both endpoints.

\Needspace{10\baselineskip}
\begin{lemma}\label{lem:straighten}
Suppose that the drawing in Proposition~\ref{prop:double} is
$k$-quasi-planar. Let $P$ and $S$ be the sphere with boundary
and its double from the proof of the proposition.
There exist a hyperbolic metric on $P$ with geodesic boundary
and a geodesic triangulation of $P$ with $8(n-2)$ faces.
Each edge arc can be replaced by a simple geodesic arc joining
the same two boundary components. The doubles of these arcs
are simple closed geodesics on $S$. Every family of distinct,
pairwise intersecting lifts of these geodesics has at most
$k-1$ members. Moreover, no triangulation vertex lies on any
of the new edge arcs.
\end{lemma}

\begin{proof}
A
\emph{pair of pants} is a sphere with three disjoint open disks
removed. Decompose $P$ into $n-2$ pairs of pants by cutting
along $n-3$ disjoint interior circles. To obtain the decomposition,
repeatedly choose two boundary circles of the remaining surface
and join them by an embedded arc whose interior avoids the
boundary. A small neighborhood of the two circles and the arc
is a pair of pants. Removing this piece decreases the number
of boundary circles of the remaining surface by one. Continue
until the final piece has three boundary circles.

Give each pair of pants a hyperbolic metric as follows. Take a
regular right-angled hyperbolic hexagon and an oppositely oriented
copy. Identify three alternate sides with their corresponding
sides in the second copy. The three pairs of unglued sides form
geodesic boundary circles: two right angles meet at each join.
The circles have equal lengths because the two hexagons are
regular and congruent. The glued sides form three arcs
perpendicular to the boundary circles.

Use the same hexagon for every pair of pants. On each boundary
circle, the two endpoints of the perpendicular arcs are separated
by half the circumference. Reassemble $P$ using orientation-reversing
isometries of matching boundary circles that identify these
endpoints. The resulting metric is hyperbolic, with geodesic
boundary. It extends smoothly across the boundary when we
double $P$ to obtain $S$.

In each hexagon, join one vertex to the three nonadjacent
vertices by geodesic segments. The three new segments divide the hexagon into
four triangles, so each pair of pants contributes eight
triangles (Figure~\ref{fig:pants}). Consequently, the triangulation
of $P$ has $8(n-2)$ faces. Every face lifts to a nondegenerate,
convex geodesic triangle.

\begin{figure}[htbp]
\centering
\begin{tikzpicture}[diagram,scale=.95]
  \draw[surface border,fill=black!5] (0,0)
    ellipse[x radius=1.65,y radius=1.2];
  \draw[surface border,fill=white] (-.65,0)
    ellipse[x radius=.35,y radius=.40];
  \draw[surface border,fill=white] (.65,0)
    ellipse[x radius=.35,y radius=.40];
  \draw[second curve] (-1.65,0)--(-1,0);
  \draw[second curve] (-.30,0)--(.30,0);
  \draw[second curve] (1,0)--(1.65,0);
  \draw[map arrow] (2.02,0)--(2.92,0);
  \foreach \h in {-.72,.72} {
    \begin{scope}[shift={(3.9,\h)}]
      \draw[surface border,fill=black!5]
        (0:.65)--(60:.65)--(120:.65)--(180:.65)--(240:.65)--(300:.65)--cycle;
      \draw[second curve] (0:.65)--(60:.65);
      \draw[second curve] (120:.65)--(180:.65);
      \draw[second curve] (240:.65)--(300:.65);
    \end{scope}
  }
  \draw[map arrow] (5.00,0)--(6.00,0);
  \begin{scope}[shift={(7.55,0)}]
    \draw[surface border,fill=black!5]
      (0:1.2)--(60:1.2)--(120:1.2)--(180:1.2)--(240:1.2)--(300:1.2)--cycle;
    \draw[second curve] (0:1.2)--(60:1.2);
    \draw[second curve] (120:1.2)--(180:1.2);
    \draw[second curve] (240:1.2)--(300:1.2);
    \foreach \a in {120,180,240}
      \draw[first curve,line width=.55pt] (0:1.2)--(\a:1.2);
  \end{scope}
  \node[align=center] at (0,-1.72) {(a) Three cutting arcs};
  \node[align=center] at (3.9,-1.72) {(b) Two hexagons};
  \node[align=center] at (7.55,-1.72) {(c) Four triangles\\in each hexagon};
\end{tikzpicture}
\caption{Three cutting arcs divide a pair of pants into two
hexagons. Each hexagon contributes four triangles. The dashed
sides are copies of the cutting arcs. The hexagons are drawn
schematically and have right angles in the hyperbolic metric.
\emph{TikZ figure prepared with OpenAI Codex assistance.}}
\label{fig:pants}
\end{figure}
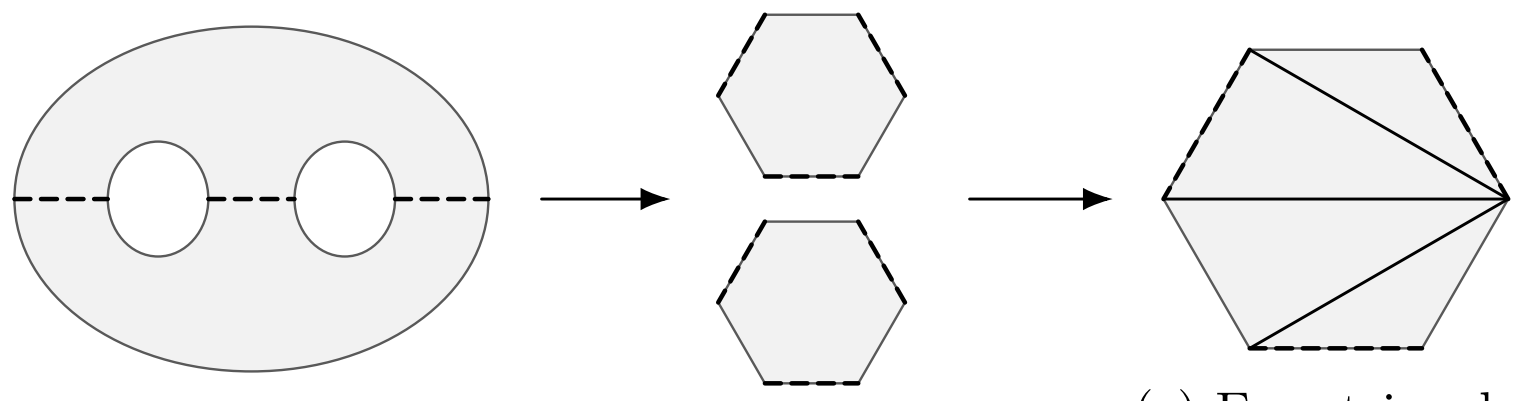

We next straighten the doubled edge curves. We use two standard
facts about closed hyperbolic surfaces: every essential simple
closed curve has a unique simple geodesic in its unoriented free
homotopy class, and distinct simple closed geodesics minimize
the number of intersections in their classes
\cite[Section~6.3]{Martelli2022}. For each edge $e$, denote by
$c_e$ the geodesic representing $\gamma_e$. Suppose that $e$ joins $u$ and $v$. The curve $\gamma_e$ has
mod-two intersection number one with $s_u$ and $s_v$.
Each gluing circle has mod-two intersection number zero with
$s_u$; for the intersection of $s_u$ with itself, take a parallel
copy. Hence $c_e$ differs from every gluing circle.

The original curve $\gamma_e$ meets $s_u$ and $s_v$ once each
and is disjoint from all other gluing circles. The geodesic
minimal-intersection property gives the same upper bounds for
$c_e$. Its odd intersection parities force at least one
intersection with each of $s_u$ and $s_v$. Thus $c_e$ meets
those two circles exactly once each and no other gluing circle.

The reflection interchanging $P_+$ and $P_-$ is an isometry that
preserves the unoriented class of $\gamma_e$. Uniqueness of the
geodesic representative implies that the reflection preserves
$c_e$. At an intersection with a gluing circle, the invariant
tangent line to $c_e$ is perpendicular to the circle. The only
other invariant line is tangent to the circle, which would force
the two geodesics to coincide. The two gluing-circle crossings
therefore cut $c_e$ into one arc in each copy of $P$, with endpoints
on $b_u$ and $b_v$. If two of the new edge arcs shared a boundary
endpoint, then their doubles would have the same tangent line
there and hence coincide, contrary to Proposition~\ref{prop:double}.

A free homotopy from $\gamma_e$ to $c_e$ pairs their lifts.
Corresponding points on the paired lifts remain at bounded
distance: each period of the lifted homotopy is a translate of
a compact fundamental strip. Consequently, the two lifts have
the same ideal endpoints.

Suppose that two geodesic lifts intersect. Their four ideal
endpoints alternate, so the corresponding lifts of the original
curves must also intersect. Indeed, these original lifts are
simple arcs joining alternating endpoints in the closed disk.
Distinct lifts of one embedded closed curve are disjoint.
Therefore, each member of a pairwise intersecting family of
geodesic lifts represents a different drawing edge. By
Proposition~\ref{prop:double}, those drawing edges cross
pairwise. The family has at most $k-1$ members.

Identify $P$ with $P_+\subset S$. Each component of the
preimage of $P_+$ in the hyperbolic plane is the intersection
of the closed half-planes on its side of the lifted boundary
circles. Hence each component is convex. A lifted edge arc
is the intersection of this convex region with its complete
geodesic.

Finally, move the triangulation vertices by a small isotopy of
$P$, and replace the displaced sides by geodesic segments with
the new endpoints. Sufficiently small moves preserve the
triangular cell decomposition. Choose all vertices off the
edge arcs.

The interior triangulation sides must also avoid crossings of edge
arcs. Work in a lift, and fix one endpoint of a side. The
other endpoint must avoid the geodesic line through the fixed
endpoint and the crossing. If the moving endpoint lies on a
boundary geodesic, that boundary meets the forbidden line in
at most one point. Otherwise the fixed endpoint would lie on
the same boundary geodesic, and the side would be a boundary
side. Thus each prohibited alignment excludes a closed set
with empty interior in the space of small vertex moves.
Only finitely many lifted crossings lie in a compact neighborhood
of any given lifted side. We may therefore avoid all such
alignments simultaneously.

Boundary sides contain no crossings of edge arcs, because
arc endpoints are distinct. Nor can a side share a segment
with a lifted edge arc: convexity would place the side's
vertices on that arc. Consequently, every intersection of a
triangulation side with an edge arc is transverse.
\end{proof}

\Needspace{6\baselineskip}
A \emph{chain} in a partially ordered set (\emph{poset}) is a
set of pairwise comparable elements. An \emph{antichain} is a
set of pairwise incomparable elements. The \emph{width} of a
finite poset is the largest cardinality of an antichain.
Two linear orders \emph{realize} a partial order if one element
precedes another in the partial order exactly when it precedes
that element in both linear orders.

Fix an orientation of $S$ and lift it to the hyperbolic plane.
For a directed side~$pq$ in the lifted triangulation, let
$\mathcal L(pq)$ consist of the complete lifts of the geodesics
$c_e$ crossing the interior of $pq$. Each member corresponds to
one of the finitely many intersections of the edge arcs with
the projected side. Distinct members may represent occurrences
of the same drawing edge at different intersections.

Every member of $\mathcal L(pq)$ is the axis of a nonidentity
covering transformation. Two distinct axes cannot share both
ideal endpoints, because the endpoints determine the geodesic.
Nor can they share exactly one ideal endpoint.

Suppose to the contrary that two axes share exactly one ideal
endpoint. In the upper half-plane model, send this endpoint to
infinity and the other endpoint of the first axis to zero.
The corresponding covering transformations then have the forms
$z\mapsto\lambda z$ and $z\mapsto\mu z+b$. We may invert the
first transformation so that $\lambda>1$; the assumption that
the axes are distinct gives $b\ne0$. Conjugating the second
transformation by the $j$th inverse power of the first yields
$z\mapsto\mu z+b\lambda^{-j}$. These transformations are
pairwise distinct and converge as $j$ increases. Their
successive quotients approach the identity, contradicting
discreteness of the covering group. Thus distinct members of
$\mathcal L(pq)$ have disjoint sets of ideal endpoints.

\begin{figure}[htbp]
\centering
\begin{tikzpicture}[diagram,scale=.92]
  \draw[surface border] (0,0) circle[radius=1.85];
  \draw[black!45,densely dotted] (-1.85,0)--(1.85,0);
  \draw[map arrow] (-1.53,0)--(1.53,0);
  \fill (-1.53,0) circle[radius=1.4pt] (1.53,0) circle[radius=1.4pt];
  \node[below=3pt] at (-1.53,0) {$p$};
  \node[below=3pt] at (1.53,0) {$q$};
  \node[left=4pt] at (-1.85,0) {$p_\infty$};
  \node[right=4pt] at (1.85,0) {$q_\infty$};
  \draw[first curve] (130:1.85)--(230:1.85);
  \draw[second curve] (100:1.85)--(310:1.85);
  \draw[black,line width=.75pt,dash pattern=on 4pt off 2pt on .8pt off 2pt]
    (60:1.85)--(260:1.85);
  \foreach \a/\lab in {130/1,100/2,60/3,230/1,310/2,260/3} {
    \fill (\a:1.85) circle[radius=1.5pt];
    \node at (\a:2.12) {$\ell_{\lab}$};
  }
  \node[align=left,anchor=west] at (3.05,.56)
    {$L_{pq}:\quad \ell_1<\ell_2<\ell_3$};
  \node[align=left,anchor=west] at (3.05,-.10)
    {$R_{pq}:\quad \ell_1<\ell_3<\ell_2$};
  \node[align=left,anchor=west] at (3.05,-.86)
    {$\ell_1\prec\ell_2,\quad\ell_1\prec\ell_3$;\\
     $\ell_2$ and $\ell_3$ are incomparable.};
\end{tikzpicture}
\caption{The endpoint orders for three geodesics in the Klein
disk. The geodesic $\ell_1$ is disjoint from both $\ell_2$ and
$\ell_3$. The latter two intersect.
\emph{TikZ figure prepared with OpenAI Codex assistance.}}
\label{fig:endpoints}
\end{figure}
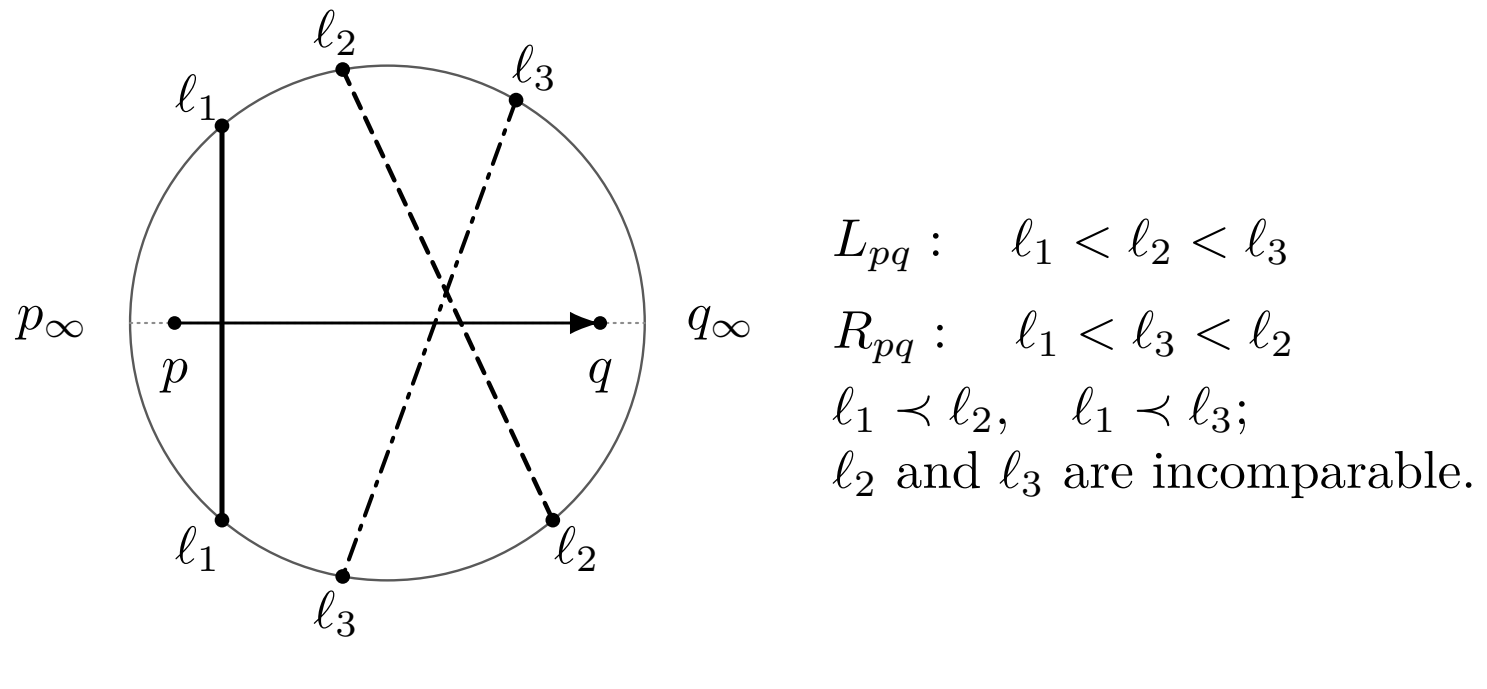

Let $p_\infty$ and $q_\infty$ denote the backward and forward
ideal endpoints of the complete line through the directed
side $pq$. They divide the ideal circle into a left arc and
a right arc, viewed from $p$ toward $q$. Every geodesic in
$\mathcal L(pq)$ has one endpoint on each arc. Read the
endpoints along each arc from $p_\infty$ to $q_\infty$.
The endpoint readings give linear orders $L_{pq}$ and $R_{pq}$ on
$\mathcal L(pq)$ (Figure~\ref{fig:endpoints}). Denote by
$Q_{pq}$ the poset realized by these orders and by $\prec$
its strict order.

For a geodesic $\ell\in\mathcal L(pq)$, let $H_q(\ell)$
denote the open half-plane bounded by $\ell$ and containing $q$.
The endpoint orders then give
\[
 \ell\prec\ell'\quad\Longleftrightarrow\quad
 H_q(\ell')\subsetneq H_q(\ell).
\]
When two geodesics crossing $pq$ are disjoint, both endpoint
orders agree with the order of their crossings along $pq$.
The later geodesic bounds the smaller half-plane containing
$q$. When two such geodesics intersect, their endpoint orders
are opposite. Consequently, two elements of $Q_{pq}$ are
incomparable exactly when their geodesics intersect.
Set $d=k-1$. Lemma~\ref{lem:straighten} implies that $Q_{pq}$
has width at most $d$. Covering transformations preserve
orientation and both endpoint orders. Thus the poset and its
two realizing orders do not depend on the choice of lift of
the directed side.

A subset~$I$ of a poset is a \emph{downset} if it contains every
element below any of its elements. It is an \emph{upset} if it
contains every element above any of its elements. In a lifted
triangle, the portions of the complete geodesics are precisely
the lifted segments of the edge arcs. We compare the orders on
the two sides crossed by such a segment.

\Needspace{7\baselineskip}
\begin{lemma}\label{lem:transport}
In a lifted triangle $pqr$, the geodesics crossing $pq$ and $pr$
form a common downset of $Q_{pq}$ and $Q_{pr}$. The geodesics
crossing $pq$ and $qr$ form a common upset of $Q_{pq}$ and
$Q_{rq}$. In both identifications, the two realizing orders agree.
\end{lemma}

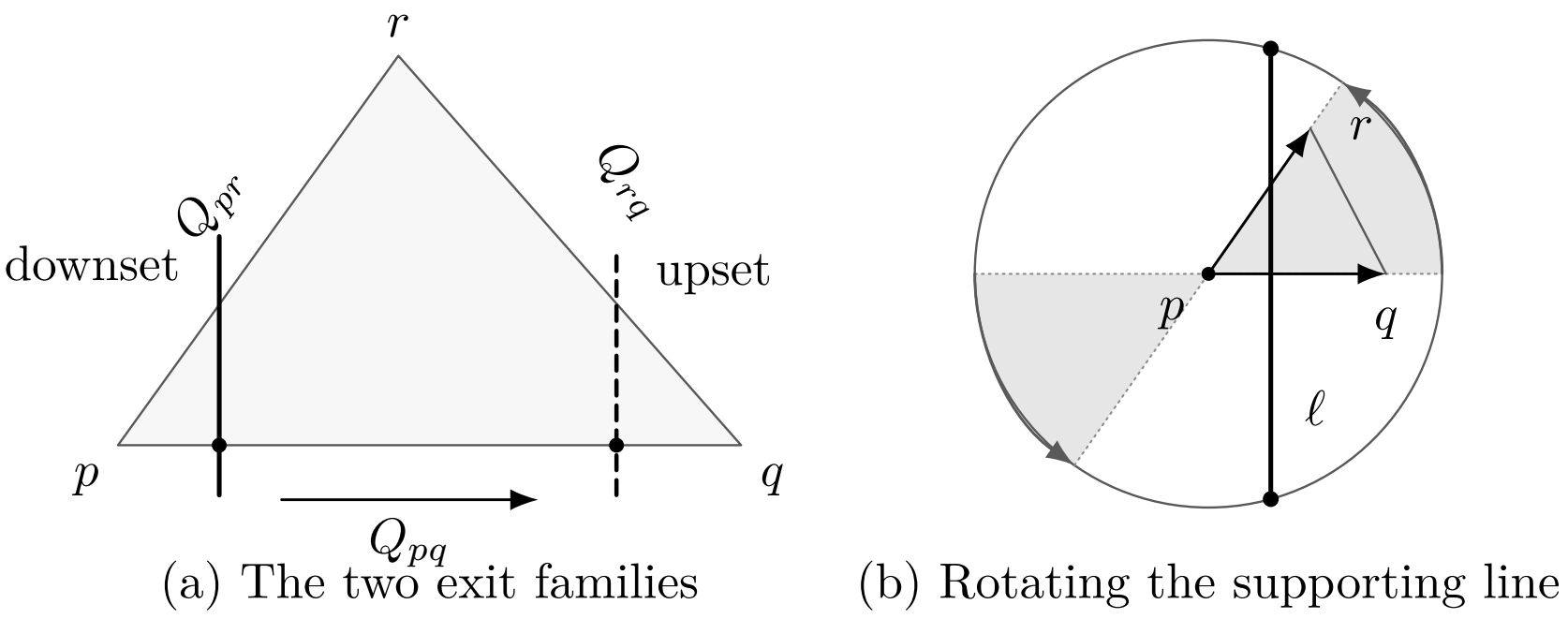
\begin{figure}[htbp]
\centering
\begin{tikzpicture}[diagram,x=1.1cm,y=1.1cm]
  \coordinate (p) at (0,0);
  \coordinate (q) at (4,0);
  \coordinate (r) at (1.8,2.5);
  \draw[surface border,fill=black!3] (p)--(q)--(r)--cycle;
  \draw[first curve] (.65,-.32)--(.65,1.34);
  \draw[second curve] (3.2,-.32)--(3.2,1.25);
  \fill (.65,0) circle[radius=1.5pt] (3.2,0) circle[radius=1.5pt];
  \node[below left] at (p) {$p$};
  \node[below right] at (q) {$q$};
  \node[above] at (r) {$r$};
  \node[left=4pt] at (.65,1.17) {downset};
  \node[right=4pt] at (3.2,1.10) {upset};
  \draw[map arrow] (1.05,-.35)--(2.70,-.35);
  \node[below] at (1.87,-.35) {$Q_{pq}$};
  \node[rotate=54] at (.60,1.55) {$Q_{pr}$};
  \node[rotate=-49] at (3.27,1.70) {$Q_{rq}$};
  \node at (2,-.90) {(a) The two exit families};
  \begin{scope}[shift={(7,1.1)}]
    \path[fill=black!10] (0,0)--(0:1.5)
      arc[start angle=0,end angle=55,radius=1.5]--cycle;
    \path[fill=black!10] (0,0)--(180:1.5)
      arc[start angle=180,end angle=235,radius=1.5]--cycle;
    \draw[surface border] (0,0) circle[radius=1.5];
    \draw[black!45,densely dotted] (180:1.5)--(0:1.5);
    \draw[black!45,densely dotted] (235:1.5)--(55:1.5);
    \draw[surface border] (0:1.14)--(55:1.14);
    \draw[map arrow] (0,0)--(0:1.14);
    \draw[map arrow] (0,0)--(55:1.14);
    \draw[first curve] (74.53:1.5)--(285.47:1.5);
    \fill (74.53:1.5) circle[radius=1.6pt]
      (285.47:1.5) circle[radius=1.6pt];
    \fill (0,0) circle[radius=1.4pt];
    \node[below left=1pt] at (0,0) {$p$};
    \node[below=3pt] at (0:1.14) {$q$};
    \node[right=4pt] at (55:1.14) {$r$};
    \node[right=3pt] at (.4,-.86) {$\ell$};
    \draw[map arrow,black!65] (0:1.5)
      arc[start angle=0,end angle=55,radius=1.5];
    \draw[map arrow,black!65] (180:1.5)
      arc[start angle=180,end angle=235,radius=1.5];
  \end{scope}
  \node at (7,-.90) {(b) Rotating the supporting line};
\end{tikzpicture}
\caption{(a) The exit through $pr$ gives a downset, and the exit
through $qr$ gives an upset. The orders on the exit sides use
directions $pr$ and $rq$. (b) As the ray $pq$ turns toward $pr$, the ideal
endpoints of its supporting line sweep the boundary arcs of the
shaded sectors. The geodesic~$\ell$ crosses both segments, and
its endpoints lie outside the swept arcs.
\emph{TikZ figure prepared with OpenAI Codex assistance.}}
\label{fig:triangle}
\end{figure}

\begin{proof}
Set $J_p=\mathcal L(pq)\cap\mathcal L(pr)$ and
$J_q=\mathcal L(pq)\cap\mathcal L(rq)$. A geodesic crossing
$pq$ separates $p$ from $q$ and avoids the vertex $r$.
It therefore crosses exactly one of the other sides of the
triangle. It belongs to $J_p$ when $r$ lies in its half-plane
containing $q$, and to $J_q$ otherwise.

If $\ell\prec\ell'$ in $Q_{pq}$ and $\ell'\in J_p$, then
$r\in H_q(\ell')\subset H_q(\ell)$, so $\ell\in J_p$.
Thus $J_p$ is a downset in $Q_{pq}$. Interchanging $q$ and $r$
gives the downset property in $Q_{pr}$. For $J_q$, if
$\ell\prec\ell'$ in $Q_{pq}$ and $\ell\in J_q$, then
$r\notin H_q(\ell)$ and hence $r\notin H_q(\ell')$, so
$\ell'\in J_q$. Interchanging $p$ and $r$ gives the upset
property in $Q_{rq}$.

To compare the realizing orders on $J_p$, rotate the ray
$pq$ toward $pr$ through the interior angle at $p$
(Figure~\ref{fig:triangle}(b)). Every intermediate ray meets the
segment $qr$. Each geodesic in $J_p$ separates $p$ from that
entire segment and therefore crosses every intermediate ray.
The supporting line of the rotating ray never shares an ideal
endpoint with a member of $J_p$. Its ideal endpoints cannot
pass those of any member of $J_p$, so the left and right
orders remain unchanged.

For $J_q$, rotate the ray $qp$ toward $qr$ about $q$. Every
geodesic in $J_q$ separates $q$ from the segment $pr$, so the
same argument preserves both orders. Reversing a directed
segment exchanges left and right and reverses each order.
Reversing both segments therefore gives agreement on $pq$
and $rq$.
\end{proof}

\section{A stable greedy partition into chains}\label{sec:greedy}

Fix a finite poset $Q$ and two specified linear orders $L$
and $R$ that realize it. Assign to each element its rank in
$R$, starting with rank one, and process the elements in order
$L$. The last element of each current chain is its \emph{tail}.
For each new element, choose the chain whose tail has the
largest rank strictly below the rank of that element, and
append the element to the chosen chain. If no such tail exists,
start a new chain. Each appended element follows its tail
in both realizing orders, so every part is a chain.

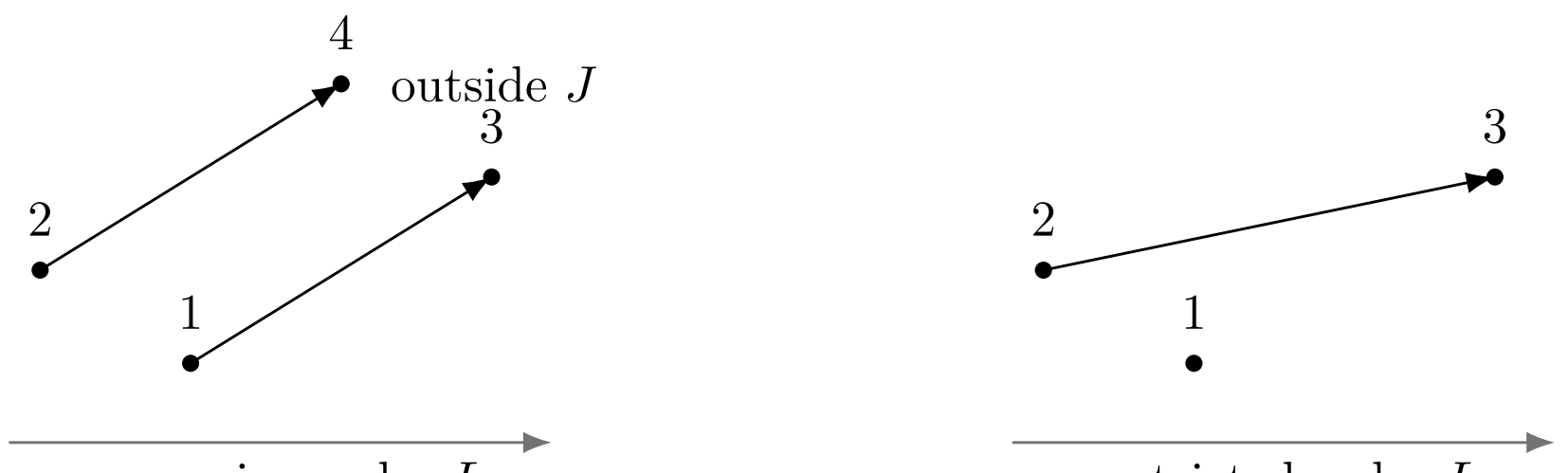
\begin{figure}[htbp]
\centering
\begin{tikzpicture}[diagram,x=1.05cm,y=.65cm]
  \begin{scope}
    \coordinate (a) at (0,2);
    \coordinate (b) at (1,1);
    \coordinate (c) at (2,4);
    \coordinate (d) at (3,3);
    \draw[map arrow] (a)--(c);
    \draw[map arrow] (b)--(d);
    \foreach \p/\lab in {a/2,b/1,c/4,d/3} {
      \fill (\p) circle[radius=1.7pt];
      \node[above=3pt] at (\p) {$\lab$};
    }
    \node[right=6pt] at (c) {outside $J$};
    \draw[map arrow,black!55] (-.2,.15)--(3.4,.15);
    \node[below] at (1.6,.15) {processing order $L$};
    \node at (1.5,-1.15) {(a) Full process};
  \end{scope}
  \begin{scope}[xshift=7cm]
    \coordinate (a) at (0,2);
    \coordinate (b) at (1,1);
    \coordinate (d) at (3,3);
    \draw[map arrow] (a)--(d);
    \foreach \p/\lab in {a/2,b/1,d/3} {
      \fill (\p) circle[radius=1.7pt];
      \node[above=3pt] at (\p) {$\lab$};
    }
    \draw[map arrow,black!55] (-.2,.15)--(3.4,.15);
    \node[below] at (1.6,.15) {restricted order $L$};
    \node at (1.5,-1.15) {(b) Process on $J$};
  \end{scope}
\end{tikzpicture}
\caption{Greedy chains for the rank sequence $2,1,4,3$. Each
arrow joins a consecutive pair in a chain. In the full process,
the element of rank $4$ replaces the tail of rank $2$. On the
downset $J=\{1,2,3\}$, the tail of rank $2$ remains available and
becomes the predecessor of the element of rank $3$.
\emph{TikZ figure prepared with OpenAI Codex assistance.}}
\label{fig:greedy}
\end{figure}

Write $\mathcal C(Q)$ for the set of ordered pairs of consecutive
elements within the resulting chains. The first element of a pair
is the \emph{predecessor} of the second, and the second is the
\emph{successor} of the first. There may be other elements
between them in the partial order. The pairs depend on the
specified realizing orders. Figure~\ref{fig:greedy} shows how
a pair can change when the rule is restricted to a downset.

\FloatBarrier

\begin{lemma}\label{lem:greedy-width}
If $Q$ has width at most $d$, then the greedy partition has at
most $d$ chains. For every downset, at most $d$ pairs in
$\mathcal C(Q)$ have their first member inside the downset and
their second member outside it.
\end{lemma}

\begin{proof}
List the current tails in decreasing order of rank. We prove
by induction on the number of processed elements that the tail
in position $j$ ends a decreasing subsequence of length $j$
in the processing order. The assertion holds when the list is
empty.

When a new tail replaces an existing one, it occupies the
same position: the replaced tail was the largest rank below
the new rank. When the new element starts a chain, its rank
is below all current tail ranks, so the new tail goes last.
For position one, the new element itself is a decreasing
subsequence of length one. At a position $j>1$, the preceding
tail has larger rank. By the induction hypothesis, it ends
a decreasing subsequence of length $j-1$; appending the new
element gives one of length $j$. Unchanged tails retain their
witnessing subsequences.

Each decreasing subsequence of ranks is an antichain, since
its elements occur in opposite orders in $L$ and $R$.
Thus there are at most $d$ tails and therefore at most $d$
chains. A chain meets a downset in an initial segment, so
at most one consecutive pair in each chain leaves the downset.
\end{proof}

\begin{lemma}\label{lem:stability}
Let $Q$ be a finite poset of width at most $d$, with two
specified realizing orders. Let $J$ be a downset or an upset,
equipped with the restrictions of those orders. Write
$\mathcal C(Q)|_J$ for the pairs in $\mathcal C(Q)$ whose two
members belong to $J$, and set $M_d=4d^2+d$.
Each of the sets
\[
 \mathcal C(Q)|_J-\mathcal C(J)
 \quad\text{and}\quad
 \mathcal C(J)-\mathcal C(Q)|_J
\]
has at most $M_d$ members.
\end{lemma}

\begin{proof}
Run the greedy process on $J$ and, separately, on $Q$.
At each stage, compare the tails of the first process with the
tails belonging to $J$ in the second process. Identify each
tail with its rank in the ambient order $R$; restricting $R$
to $J$ preserves all rank comparisons. For two sets $A$ and
$B$ of tail ranks, put
\[
 U=A\cup B\qquad\text{and}\qquad D=(A-B)\cup(B-A).
\]
For a rank $t\in U$, let $\rho_U(t)=|\{u\in U:u\le t\}|$
denote its position in the increasing list of $U$.
Define the nonnegative potential
\begin{equation}\label{eq:potential}
 \Phi(A,B)=\sum_{t\in D}\rho_U(t).
\end{equation}
A tail is \emph{common} when its rank belongs to both sets,
and \emph{exclusive} when its rank belongs to exactly one.
A \emph{common update} processes the same new rank $v$ in both
sets. In each set, the greedy rule replaces the largest tail
below $v$ by $v$, or adds $v$ if no such tail exists.

\Needspace{9\baselineskip}
The potential $\Phi$ satisfies the following three bounds
\begin{enumerate}
\item If the two predecessors exist and agree, a common update
does not increase $\Phi$. If neither predecessor exists, the
update increases $\Phi$ by exactly $|D|$.
\item If the predecessors differ, including the case when only
one exists, a common update decreases $\Phi$ by at least one.
\item Adding or deleting one tail in one set increases $\Phi$
by at most the larger of the sizes of the union before and
after the change.
\end{enumerate}
For the first bound, replacing a common predecessor by a
larger common tail leaves $D$ unchanged. No element of $D$
can move to a higher position in $U$. If neither predecessor
exists, the newly added common tail lies below every old tail.
It increases the position of each exclusive tail by one,
giving an increase of exactly $|D|$.

For the second bound, first suppose that both predecessors
exist and are different. Denote them by $a$ and $b$, choosing
the labels so that $a>b$. The tail $a$ is exclusive: if the
other set contained $a$, it would choose $a$ instead of $b$.
There are two cases. If $b$ is common, the update removes $a$
and makes $b$ exclusive. The union replaces $a$ by $v>a$.
The position of $b$ remains fixed and is smaller than the
former position of $a$. No other exclusive tail moves to a
higher position. If $b$ is exclusive, both $a$ and $b$
disappear, and $v$ becomes common. Replacing two smaller
ranks by a larger one cannot increase the contribution of
any remaining exclusive tail.

Now suppose that only one predecessor exists. It must be
exclusive, and the update removes it. In the union, the removed
rank is replaced by the larger common rank $v$. Again, no
remaining exclusive tail gains a higher position. In all three
cases, the potential decreases by at least one.

For the third bound, consider first the insertion of a tail.
If its rank already belongs to the other set, the tail becomes
common, and its contribution disappears. Otherwise, inserting
a new rank $t$ adds its position in the new union and raises
by one the position of each exclusive rank above it. The total
increase is at most the size of the new union. Next consider
a deletion. Deleting a common tail from one set makes it
exclusive, adding its position without changing the union.
Deleting an exclusive tail removes its contribution and cannot
increase any other position.

The induced poset on $J$ has width at most $d$. Both tail
sets start empty and have at most $d$ elements by
Lemma~\ref{lem:greedy-width}. The full process also handles elements outside $J$.

\emph{Case 1: $J$ is a downset.}
Every predecessor of an element of $J$ belongs to $J$.
Thus, when the full process handles an element of $J$, its
tails inside $J$ undergo the ordinary greedy update.
Processing an element outside $J$ can remove a tail from $J$
only when a chain leaves the downset. Each of the at most
$d$ chains can leave $J$ at most once. Hence at most $d$
such deletions occur.

\emph{Case 2: $J$ is an upset.}
Processing an element outside $J$ cannot affect a tail in $J$,
because every element extending such a tail also belongs to
$J$. When the full process handles an element of $J$, the
tails inside $J$ undergo an ordinary greedy update unless
the full process chooses a predecessor outside $J$. This
exception marks the first entry of a full chain into $J$.
Consequently, there are at most $d$ exceptions.

At an exceptional entry, the full process adds the new tail
inside $J$ without removing any tail there. To simulate this,
perform the ordinary greedy update and reinsert its removed
tail, if one exists. The intermediate tail set is no larger
than the actual tail set after the update, so its size is
also at most $d$.

In either case, there are at most $d$ extra insertions or
deletions between common updates. By the third bound,
these changes increase $\Phi$ by at most $2d^2$ in total. At most
$d$ common updates have no predecessor in either set: each
such update starts a chain in the ordinary process on $J$.
By the first bound, each increases $\Phi$ by at most $2d$.
Since $\Phi$ begins at zero and remains nonnegative, the
second bound limits the number of common updates with
different predecessors by $4d^2$. Here a missing predecessor
counts as a difference.

At an exceptional entry into an upset, the simulated update
may select a predecessor in $J$ even though the full process
selects one outside $J$. Returning to the actual process can
therefore introduce at most one additional discrepancy per
exception. There are at most $d$ exceptions. Thus at most
$4d^2+d$ elements of $J$ have different predecessors within
$J$ in the two processes. Every consecutive pair lost on
restriction has a distinct second element among these, and
the same holds for every gained pair.
\end{proof}

\begin{lemma}\label{lem:change}
Suppose that $Q_1$ and $Q_2$ have width at most $d$ and share a
downset $J$, or share an upset $J$. Put $B_d=8d^2+2d$.
If their two realizing orders agree on $J$, then at most $B_d$
pairs of $\mathcal C(Q_1)$ inside $J$ are absent from
$\mathcal C(Q_2)$.
\end{lemma}

\begin{proof}
Because the two realizing orders agree on $J$, their restrictions
give the same greedy partition of $J$. Consider a pair in
$\mathcal C(Q_1)|_J$ that is absent from $\mathcal C(Q_2)|_J$.
Either the pair is absent from $\mathcal C(J)$, or it belongs
to $\mathcal C(J)-\mathcal C(Q_2)|_J$.
Lemma~\ref{lem:stability} bounds each set of pairs by $M_d$.
There are therefore at most $2M_d=B_d$ such pairs.
\end{proof}

\section{Counting the edges}\label{sec:count}

\begin{proof}[Proof of Theorem~\ref{thm:main}]
Let $G$ be a graph with the drawing specified in
Theorem~\ref{thm:main}, and put $d=k-1$. Use the geodesic
arcs and the triangulation from Lemma~\ref{lem:straighten}.
On each directed side $pq$, run the greedy rule in the
processing order $L_{pq}$, using ranks from $R_{pq}$.
We count edges by recording chain ends and consecutive pairs that cease to be consecutive.

\Needspace{13\baselineskip}
Fix a triangle and a direction on one of its sides. Choose
a lift of this triangle, denote its directed side by $pq$,
and let $r$ be its third vertex. An element of $Q_{pq}$
represents an occurrence of a drawing edge crossing that side.
A \emph{record} consists of the triangle, the directed side,
and an occurrence satisfying one of the following conditions:
\begin{enumerate}
\item The occurrence is the last element of its greedy chain.
\item The occurrence has a successor, but the corresponding
two arcs leave the triangle $pqr$ through different sides.
\item The occurrence has a successor, and the corresponding
two arcs leave through the same side, but their ordered pair
is absent from the greedy partition on the outgoing directed side.
\end{enumerate}
For a common exit, use direction $pr$ or $rq$, according to
the exit side. There are at most $d$ records of Type~1,
one at the end of each chain. For Type~2, a consecutive pair
with different exits leaves the downset of arcs exiting
through $pr$. Lemma~\ref{lem:greedy-width} allows at most
$d$ such pairs. For Type~3, Lemma~\ref{lem:transport}
identifies each exit family with a common downset or a common
upset in the incoming and outgoing posets.
Lemma~\ref{lem:change} allows at most $B_d$ lost pairs
in each family. Since there are two exit families, a single
directed side contributes at most $2d+2B_d$ records.

We show that every drawing edge labels at least one record.
Fix an edge $e=uv$ and start at the unique point where its
geodesic arc meets $b_u$. The point lies in the interior
of a boundary side, since no triangulation vertex is on an
edge arc. Choose either direction on that side. If the
occurrence of $e$ is last in its chain, it produces a
Type~1 record. Otherwise the successor represents a different
edge $f$ incident with $u$: each such edge meets $b_u$
only once. Let $w$ be the other endpoint of $f$. Since
$G$ is simple, the vertices $w$ and $v$ are distinct.

Follow the arcs of $e$ and $f$ through the adjacent triangle.
If they leave through different sides, the occurrence of
$e$ produces a Type~2 record. They cannot leave together
through a boundary side, since that would put their other
endpoints on the same boundary circle, contrary to $w\ne v$.

Suppose that both arcs leave through the same interior side.
Compare their ordered pair with the greedy partition on that
outgoing directed side. If the pair is absent, the occurrence
of $e$ produces a Type~3 record. If it is present, the same
consecutive pair enters the next triangle, with $e$ still the
first element. Retain the direction on the shared side when
entering the next triangle. The greedy partition does not
depend on the chosen lift of that side.

Each arc crosses only finitely many triangulation sides.
If no record were encountered, the two arcs would continue
together to a common boundary exit. As shown above, this is
impossible. Thus every edge labels at least one record.

\Needspace{7\baselineskip}
Every record has one edge label, and every edge labels at
least one record. Each triangle has six directed sides, and
the triangulation has $8(n-2)$ triangles. Since $d\ge 1$,
the record count gives
\begin{align*}
 |E(G)|&\le 48(n-2)(2d+2B_d)\\
       &= (768d^2+288d)(n-2)\\
       &\le 1056d^2(n-2).\qedhere
\end{align*}
\end{proof}

\Needspace{22\baselineskip}
\section*{Funding}

The author received no funding for this work.

\section*{Declaration of competing interests}

The author declares no competing interests.

\section*{Declaration of generative AI and AI-assisted technologies in the manuscript preparation process}
GPT-5.6 Sol (OpenAI) identified the
paper by Beeker and Lazarovich and suggested the connection to
quasi-planar graphs. The author studied the connection, wrote
the original proof of fixed-$k$ linearity by hand, and used
GPT-5.6 Sol to transcribe photographs of that proof into LaTeX.
GPT-6 later proposed the doubling construction and the quantitative
argument, including the greedy chain
partition and stability estimate. The author used Codex to prepare the TikZ figures.
The author is responsible for the mathematical content and proofs.

\Needspace{9\baselineskip}

\end{document}